\documentclass[12pt]{amsart}
\usepackage{amssymb,latexsym,comment,url}
\usepackage{cellspace} 
\newcommand{\rk}{\mathrm{rank}}

\newcommand{\Q}{{\mathbb Q}}

\newenvironment{Proof}{\par\noindent{\sc Proof:}}%
                      {\hspace*{\fill}\nobreak$\Box$\par\medskip}
                       {\hspace*{\fill}\nobreak$\Box$\par\medskip}

\newtheorem{Proposition}{Proposition}[section]
\newtheorem{Theorem}[Proposition]{Theorem}

\theoremstyle{definition}
\newtheorem{Definition}[Proposition]{Definition}
\newtheorem{Remark}[Proposition]{Remark}
 
\renewcommand{\baselinestretch}{1.1}

\begin{document}
\title{Strong rational Diophantine $D(q)$-triples}
\author[A. Dujella]{Andrej Dujella}
\address{Department of Mathematics, Faculty of Science, University of Zagreb, Bijeni\v{c}ka
cesta 30, 10000 Zagreb, Croatia}
\email{duje@math.hr}
\author[M. Paganin]{Matteo Paganin}
\address{Foundation Development Program,
 Sabanc{\i} University,
  Tuzla, \.{I}stanbul, 34956 Turkey}
\email{matteo.paganin@sabanciuniv.edu}
\author[M. Sadek]{Mohammad Sadek}
\address{Faculty of Engineering and Natural Sciences,
 Sabanc{\i} University,
  Tuzla, \.{I}stanbul, 34956 Turkey}
\email{mohammad.sadek@sabanciuniv.edu}

\begin{abstract}
We show that for infinitely many square-free integers $q$ there exist infinitely many triples of rational
numbers $\{a,b,c\}$ such that $a^2+q$, $b^2+q$, $c^2+q$, $ab+q$, $ac+q$ and $bc+q$ are squares of rational
numbers.
\end{abstract}

\subjclass[2010]{11D09, 11G05}
\keywords{strong Diophantine triples, elliptic curves.}

\maketitle

\section{Introduction}

Classically, a \emph{Diophantine $m$-tuple} is a set $\{a_1,\ldots,a_m\}$ of $m$ non-zero integers with the property that $a_ia_j+1$ is a square, whenever $i\ne j$; such an $m$-tuple is called \emph{rational} if we allow its elements to be non-zero rational numbers.

Fermat found the first Diophantine quadruple in integers $\{1,3,8,120\}$.
In 1969, Baker and Davenport \cite{BD}
proved that Fermat's set cannot be extended to a Diophantine quintuple.
This result motivated the conjecture that there does not exist a Diophantine quintuples in integers.
The conjecture has been proved recently by He, Togb\'e and Ziegler \cite{HTZ}.

The first example of a rational Diophantine quadruple, the set
$\{\frac{1}{16},\, \frac{33}{16},\, \frac{17}{4},\, \frac{105}{16}\}$ was found by Diophantus, while
Euler proved that there exist
infinitely many rational Diophantine quintuples (see \cite{Hea}),
In 1999, Gibbs found the first example of rational Diophantine sextuple
$\{ \frac{11}{192}, \frac{35}{192}, \frac{155}{27}, \frac{512}{27}, \frac{1235}{48}, \frac{180873}{16} \}$
(see \cite{Gibbs1}).
In 2017, Dujella, Kazalicki, Miki\'c and Szikszai \cite{DKMS} proved that there are
infinitely many rational Diophantine sextuples, and alternative constructions of families
of rational Diophantine sextuples are given in \cite{Duje-Matija}, \cite{DKP-sext} and \cite{DKP-split}.
It is not known whether there exist any rational Diophantine septuple.
More information on Diophantine $m$-tuples can be found in the survey article \cite{Duje-notices}.

Dujella and Petri\v{c}evi\'c in \cite{DP} introduced the notion of \emph{strong} rational Diophantine $m$-tuple, as a rational Diophantine $m$-tuple with the additional property that $a_i^2+1$ is a rational square for every $i=1,\ldots,m$. They proved that there exist infinitely many strong rational Diophantine triples.
One such example is the set $\{1976/5607, 3780/1691, 14596/1197\}$.

Let $q$ be a rational number. A set $\{a_1, \ldots, a_m\}$ of nonzero integers (rationals)
is called a (rational) $D(q)$-$m$-tuple, if $a_ia_j + q$ is a square of a rational number for all
$1 \leq i < j \leq m$. It is known that for every rational number $q$
there exist infinitely many rational $D(q)$-quadruples, and that there are infinitely many
square-free integers $q$ for which there exist infinitely many rational $D(q)$-quintuples
(see \cite{D-quint,DF}).

In this paper, we will consider the problem which arises if we combine the two above mentioned
variants of Diophantine $m$-tuples.

\begin{Definition}
Let $q$ be a rational number. A \emph{strong rational Diophantine $D(q)$-$m$-tuple}
is a set of non-zero rationals $\{a_1,\ldots,a_m\}$ such that $a_ia_j+q$ is a square for
all $i,j=1,\ldots,m$ (including the case $i=j$).
\end{Definition}

As we already mentioned, the case $q=1$ was studied in \cite{DP}.
The case $q=-1$ was studied in \cite{DGPT} and it was shown that there exist
infinitely many strong rational $D(-1)$-triples (in \cite{DGPT} they are called
strong Eulerian triples because of the direct connection between $D(-1)$-$m$-tuples
and so called Eulerian $m$-tuples, which are sets with property that $xy+x+y=(x+1)(y+1)-1$ is a perfect
square for all elements $x,y$ of the set).

Our main result is the following theorem.

\begin{Theorem}\label{Theorem}
There exist infinitely many square-free integers $q$ with the property
that there exist infinitely many strong rational Diophantine $D(q)$-triples.
\end{Theorem}


\section{Construction of strong rational Diophantine $D(q)$-pairs and triples}


One may see easily that if $\{a_1,\ldots,a_m\}$ is a strong rational Diophantine $D(q)$-$m$-tuple, then $\{za_1,\ldots,za_m\}$ is a strong rational Diophantine $D(z^2q)$-$m$-tuple.
Therefore, it is enough to consider the problem of existence of strong rational
Diophantine $D(q)$-triples for square-free integers $q$, and will do so in Section \ref{sec:3}.
Also, since we may choose $z=1/a_1$ there is no lost of generality if we assume that $a_1=1$
and consequently $1^2 + q = r^2$, i.e.  $q=r^2-1$.


We will now explain a construction of strong rational Diophantine $D(q)$-pairs which
use properties of related elliptic curves.

\begin{Proposition}\label{Proposition2}
For all rational numbers $r$, $r\neq 0, \pm 1, \pm \frac{1}{2}$,
there exist infinitely many rational numbers $b$ such that $\{1,b\}$ is a strong rational Diophantine $D(r^2-1)$-pair.
\end{Proposition}
\begin{Proof}
%
%
%
For convenience, we set $q=r^2-1$. We consider the elliptic curve $E^q$ defined by the equation
\[
E^{q}:\quad y^2 = (x+q)(x^2+q) = x^3+qx^2+qx+q^2.
\]
The curve $E^q$ is non-singular for $q\neq 0, -1$, i.e. for $r\neq 0, \pm 1$,
so in what follows we will always assume that $r\neq 0, \pm 1$.
Some obvious rational points on $E^{q}$ are
\[
T^{q} = (-q,0),\qquad P^{q} = (0,q),\qquad S^{q} = (1,1+q).
\]
It is easily checked that $T^{q}+P^{q}+S^{q}=O$.

Any rational number $b$ such that $\{1,b\}$ is a strong rational Diophantine $D(q)$-pair, is the $x$-coordinate of a point on $E^{q}$.

Standard $2$-descent (see e.g. \cite[4.2, p.85]{Knapp})
yields that the $x$-coordinate $b$ of any point in $2E^{q}(\Q)$ satisfies that $\{1,b\}$ is a strong rational Diophantine $D(q)$-pair. Hence, we will finish the proof
if we show that $E^{q}$ has rank at least 1 over $\mathbb Q$.

We notice that
\[
2S^{q} = \left(\frac14+q,\frac{q}2+\frac18\right) = \left(\frac54-r^2,\frac{r^2}2-\frac38\right).
\]
Assume for the moment that $r$ is an integer.
Since the $y$-coordinate of $2S^{q}$ cannot be an integer, by
the Lutz-Nagell theorem $S^{q}$ has infinite order and $\rk(E^{q})$ is at least 1.
Let us consider now the general case when $r$ is a rational number.
We want to show that again the point $S^{q}$ has infinite order.
By Mazur's classification of torsion points of elliptic curves over $\mathbb Q$, it is enough to check that $kS^{q}$ is not the point at infinity for
$k\leq 12$ by considering rational roots of the denominators of the coordinates.
We obtain that the only rational roots of denominators are $r=\pm \frac{1}{2}$, in which cases
the point $S^{q}$ is of order $3$. For all other rational numbers $r$, the point $S^{q}$ is of infinite order.
\end{Proof}

 By the proof of Proposition~\ref{Proposition2},
 we may use the $x$-coordinate of $2kS^{r^2-1}$, $k$ is an integer, to construct families of strong rational Diophantine $D(r^2-1)$-pairs. However, since the $x$-coordinate of $S^{r^2-1}$ (which is equal to $1$)
 satisfies that conditions that both $x+r^2-1$ and $x^2+r^2-1$ are rational squares, by $2$-descent, we conclude that we may also
 use the  $x$-coordinate of $(2k+1)S^{r^2-1}$.

 For example, the $x$-coordinates of $2S^{r^2-1}$, $3S^{r^2-1}$ and $4S^{r^2-1}$
 yield that the pairs
\[
\left\{1, \frac54-r^2\right\}, \qquad
\left\{1, \frac{-16r^4+16r^2+1}{16r^4-8r^2+1}\right\}, \qquad \left\{1, \frac{256r^8-768r^6+864r^4-496r^2+145}{256r^4-384r^2+144}\right\}
\]
are $D(r^2-1)$-pairs.

By extending the first of these three families of pairs, we will construct infinitely many strong rational Diophantine $D(r^2-1)$-triples for rational numbers $r$ of certain form.
More precisely, we prove the following proposition.

\begin{Proposition}\label{Proposition3}
For any rational number $t$ different from $0$, $\displaystyle\pm \frac15$, $\displaystyle\pm \frac35$, $\displaystyle\pm \frac75$ or $\displaystyle\pm\frac7{15}$, the triple
\[
\left\{1, -\frac{625t^4-930t^2+49}{1024\,t^2}, -\frac{(5t+1)(5t-1)(5t+7)(5t-7)}{1600\,t^2}\right\}
\]
is a strong rational Diophantine $D(q)$-triple, with
\[
q=\frac{(t-1)(t+1)(25t+7)(25t-7)}{1024\,t^2}.
\]
\end{Proposition}
\begin{Proof}
In what follows we will use the symbol $\square$ to mean a square of a rational number.
A strong rational Diophantine $D(q)$-triple $\{a,b,c\}$ amounts to the following conditions being simultaneously verified:
\[
\begin{array}{lll}
a^2+q=\square_{aa}, & b^2+q=\square_{bb}, & c^2+q=\square_{cc}, \\
ac+q=\square_{ac}, & ab+q=\square_{ab}, & bc+q=\square_{bc}.
\end{array}
\]
We set $q=r^2-1$, $a=1$, and $\displaystyle b=\frac54-r^2$, for a rational number
$r\ne 0, \pm1, \pm \frac{1}{2}$.

We want to find $c$, different from 1 and $b$, such that $\{1,b,c\}$ is a strong Diophantine $D(q)$-triple. From the condition $c+q=s^2$, we shall write $c=s^2-r^2+1$, for some rational number $s$.
We search for such values of the form $s=kr$. The condition $bc+q=\square_{bc}$ then becomes
\[
p(k,r)=\frac54r^2k^2-r^4k^2-\frac54r^2+r^4+\frac14=\square_{bc}.
\]
This is possible for the values of $k$ that make the discriminant of $p(k,r)$ vanish. The discriminant of $p(k,r)$, with respect to $r$, is equal to
\[
-\frac1{64}(5k-3)^2(5k+3)^2(k-1)^3(k+1)^3,
\]
so to have $c\neq 1$ we can choose $k=3/5$.
Then $p(3/5,r)=\left(\frac{8r^2-5}{10}\right)^2$. Thus, the only condition left
is $c^2+q=\square_{cc}$, with $\displaystyle c=-\frac{16}{25}r^2+1$, that translates into
\[
\frac1{625}r^2(256r^2-175) = \square_{cc}.
\]
This implies that we need to find $t\in\Q$ such that
\[
(256r^2-175) = (16r + 25t)^2,
\]
that results into the equality $\displaystyle r=-\frac1{32}\frac{25t^2+7}t$. Substituting this value in the formulas for $b$, $c$, and $q$, we finally obtain that the triple
\[
\left\{1, -\frac{625t^4-930t^2+49}{1024\,t^2}, -\frac{(5t+1)(5t-1)(5t+7)(5t-7)}{1600\,t^2}\right\}
\]
is a strong rational Diophantine $\displaystyle D\left(\frac{(t-1)(t+1)(25t+7)(25t-7)}{1024\,t^2}\right)$-tuple. Finally, the two elements different from 1 are distinct if and only if $t$ is different from $\displaystyle\pm \frac35$ or $\displaystyle\pm\frac7{15}$.\end{Proof}

\section{Proof of Theorem \ref{Theorem}} \label{sec:3}


From Proposition~\ref{Proposition3}, we only need to prove that, for infinitely many square-free integers $q$, there are infinitely many rational numbers $t$ such that
\[
\frac{(t-1)(t+1)(25t+7)(25t-7)}{1024\,t^2} = q w^2,
\]
for some rational number $w$.
Then by dividing all elements of the triple from Proposition \ref{Proposition3} by $w$,
we get a strong rational $D(q)$-triple.

In other words, we need to study the quartic curve $Q_q: qv^2=(t-1)(t+1)(25t+7)(25t-7)$. The latter curve is the $n$-quadratic twist of the curve $Q:v^2=(t-1)(t+1)(25t+7)(25t-7) $.

The quartic curve $Q$ is birationally equivalent, by the substitutions
 $t=\frac{144+x}{144-x}$, $v=\frac{576 y}{(x-144)^2}$,
 to the elliptic curve described by the Weierstrass equation:
\[
E_1: y^2 = x(x+81)(x+256);
\]
similarly, $Q_q$ is birationally equivalent, by the same substitutions, to 
$$ E_q: qy^2 = x(x+81)(x+256). $$

We will conclude our proof if we can find infinitely many square-free $q$ for which $\rk(E_q)\ge 1$.
We will follow the reasoning from \cite{DF}.
It is well-known (see e.g. \cite{DGL}) that for the elliptic curve given by the equation
$y^2 = f(x)$, the point $(u,1)$ is a rational point of infinite order in
$E_{f(u)}(\Q)$. By writing $u=u_1/u_2$, we get that
for all integers $q$ of the form
\begin{equation} \label{eq:q}
q=u_1u_2(u_1+81u_2)(u_1+256u_2)
\end{equation}
the curve $E_q$ has positive rank.
This gives us infinitely many square-free values of $q$ for which the rank is positive,
and thus for which there exist infinitely many strong rational $D(q)$-quintuples.
Indeed, for fixed $\varepsilon > 0$ and sufficiently large
$N$, there are at least $N^{1/2-\varepsilon}$
square-free integers $q$, $|q| \leq N$, of the form (9) (see e.g. [GM]).

\qed

\section{Examples and remarks}

We computed the rank of $E_q$ for small values $q$ by {\tt mwrank} \cite{Cremona},
and obtained that rank is positive for the following square-free integers in the range $|q|<100$:
\begin{gather*}
 \mbox{-5}, \mbox{-6}, \mbox{-7}, \mbox{-11}, \mbox{-17}, \mbox{-19}, \mbox{-21}, \mbox{-22}, \mbox{-23},
\mbox{-29}, \mbox{-30}, \mbox{-34}, \mbox{-35}, \mbox{-37}, \mbox{-38}, \mbox{-39}, \mbox{-43}, \mbox{-46},
\mbox{-51}, \mbox{-55}, \\
 \mbox{-57}, \mbox{-58}, \mbox{-61}, \mbox{-62}, \mbox{-66}, \mbox{-67}, \mbox{-69},
\mbox{-74}, \mbox{-77}, \mbox{-78}, \mbox{-79}, \mbox{-83}, \mbox{-85}, \mbox{-86}, \mbox{-87}, \mbox{-91},
\mbox{-93}, \mbox{-94}, \mbox{-95}, \\
 2, 6, 10, 13, 15, 17, 23, 26, 29, 30, 31, 33, 35, 37, 42, 46\\
 47, 53, 55, 58, 59, 66, 69, 77, 78, 79, 82, 91, 93, 95.
\end{gather*}

In next table we give some examples of strong rational $D(q)$-triples $\{a,b,c\}$,
for small values of $q$,
obtained by the construction from Theorem \ref{Theorem}. We provide also the corresponding
parameter $t$. \medskip

\begin{center}
\begin{tabular}{|c|c||Sc|Sc|Sc|}
\hline
$t$ & $q$ & $a$ & $b$ & $c$ \\
%
\hline
$\displaystyle\frac{37}{125}$ & $-11$ & $\displaystyle\frac{370}{27}$ & $\displaystyle\frac{21122}{4995}$ & $\displaystyle\frac{75578}{13875}$ \\
\hline
$\displaystyle\frac{11}{25}$ & $-7$ & $\displaystyle\frac{44}{9}$ & $\displaystyle\frac{1051}{396}$ & $\displaystyle\frac{736}{275}$ \\
\hline
$\displaystyle\frac{101}{155}$ & $-6$ & $\displaystyle\frac{3131}{684}$ & $\displaystyle\frac{21031705}{8566416}$ & $\displaystyle\frac{591745}{237956}$ \\
\hline
$\displaystyle-\frac{23}{25}$ & $-5$ & $\displaystyle\frac{23}{3}$ & $\displaystyle\frac{709}{276}$ & $\displaystyle\frac{1827}{575}$ \\
\hline
$\displaystyle-\frac{119}{457}$ & 2 & $\displaystyle\frac{7769}{1638}$ & $\displaystyle\frac{38893009}{50902488}$ & $\displaystyle\frac{50817649}{35348950}$ \\
\hline
$\displaystyle-\frac{23}{265}$ & $6$ & $\displaystyle-\frac{1219}{1188}$ & $\displaystyle\frac{32386295}{5792688}$ & $\displaystyle\frac{542735}{160908}$  \\
\hline
$\displaystyle\frac{1}{31}$ & $10$ & $\displaystyle\frac{31}{66}$ & $\displaystyle-\frac{173279}{8184}$ & $\displaystyle-\frac{229437}{17050}$ \\
\hline
$\displaystyle\frac{1}{25}$ & $13$ & $\displaystyle\frac{2}{3}$ & $\displaystyle-\frac{58}{3}$ & $\displaystyle-\frac{306}{25}$ \\
\hline
\end{tabular}
\end{center}\bigskip

Just for fun, we also give a triple for $q=2019$:
{\tiny
\begin{eqnarray*}
a\hspace{-5pt} &=\hspace{-5pt}& -\frac{108425648984099462722723028577175690286281358594075905}{1979956008273178460383709106649735645388794922519592},\\
b\hspace{-5pt} &=\hspace{-5pt}& \frac{19903622160350297465727113805280431196879309182571712631429120369343905672609842407986879203598345282474239}{858712060627945518172033052697448822731672169127032763561281839945494931723647684264003999284669990523040}, \\
c\hspace{-5pt} &=\hspace{-5pt}& \frac{2314875761476160622113200620592571545156501721172189311604105086986000693279887159122625184996952958005759}{596327819880517720952800731039895015785883450782661641362001277739927035919199780738891666169909715641000}.
\end{eqnarray*}
}

\begin{Remark}
In Theorem \ref{Theorem}, the existence of infinitely many square-free integers $n$ for which there are infinitely many $D(n)$-triples mounts down to investigating the Mordell-Weil rank of the quadratic twists 
$E_q: qy^2 = x(x+81)(x+256)$. Goldfeld's minimalist conjecture asserts that for $50\%$ of square-free integers $q$, one would expect that $\rk (E_q)$ is positive, hence there are infinitely many strong rational
$D(q)$-triples for at least $50\%$ of square-free integers $q$.
See \cite{DF} for reasoning how the Parity Conjecture
implies that for $q$'s in certain arithmetic progressions the rank of $E_q$ is odd, and hence positive.
\end{Remark}

\begin{Remark}
Note that the elliptic curve $E_1$, given by the equation $y^2 = x(x+81)(x+256)$
has rank $0$ and torsion group $\mathbb{Z}/2\mathbb{Z} \times \mathbb{Z}/8\mathbb{Z}$.
For curves with such torsion group it is known that there are infinitely many quadratic twists
with rank $\geq 4$ (see \cite{M,RS}).
\end{Remark}

{\bf Acknowledgements:}
A. D. acknowledges support from the QuantiXLie Center of Excellence, a project co-financed
by the Croatian Government and European Union through the European Regional Development
Fund - the Competitiveness and Cohesion Operational Programme (Grant KK.01.1.1.01.0004). A.
D. was also supported by the Croatian Science Foundation under the project no. IP-2018-01-1313.
This work started while A. D. visited Faculty of Engineering and Natural Sciences,
 Sabanc{\i} University in November 2019.
 He thanks the colleagues of this institution for their hospitality and support. M. S. was supported by the starting project B.A.CF-19-01964, Sabanc{\i} University.


\begin{thebibliography}{99}
\frenchspacing
\renewcommand{\baselinestretch}{1}

\bibitem{BD}
A. Baker and H. Davenport, \textit{The equations $3x^2 - 2 = y^2$ and $8x^2 - 7 = z^2$},
Quart. J. Math. Oxford Ser. (2) {\bf 20} (1969), 129--137.

\bibitem{Cremona}
J. Cremona, Algorithms for Modular Elliptic Curves, Cambridge University Press, Cambridge, 1997.

\bibitem{D-quint}
A. Dujella, \textit{A note on Diophantine quintuples}, in: Algebraic Number Theory and Diophantine Analysis
(F. Halter-Koch, R. F. Tichy, eds.), Walter de Gruyter, Berlin, 2000, pp. 123--127.


\bibitem{Duje-notices}
A. Dujella, \textit{What is ... a Diophantine $m$-tuple?}, Notices Amer. Math. Soc. {\bf 63} (2016), 772--774.

\bibitem{DF}
A. Dujella and C. Fuchs, \textit{On a problem of Diophantus for rationals},
J. Number Theory {\bf 132} (2012), 2075--2083.

\bibitem{DGL}
A. Dujella, I. Gusi\'c and L. Lasi\'c, \textit{On quadratic twists of elliptic curves
$y^2 = x(x - 1)(x - \lambda)$}, Rad Hrvat. Akad. Znan. Umjet. Mat. Znan. {\bf 18} (2014), 27--34.

\bibitem{DGPT}
A. Dujella, I. Gusi\'c, V. Petri\v cevi\'c and P. Tadi\'c,
\textit{Strong Eulerian triples},
Glas. Mat. Ser. III {\bf 53} (2018), 33-42.

\bibitem{Duje-Matija}
A. Dujella and M. Kazalicki, \textit{More on Diophantine sextuples}, in: Number Theory - Diophantine problems, uniform distribution and applications, Festschrift in honour of Robert F. Tichy's 60th birthday (C. Elsholtz, P. Grabner, Eds.), Springer-Verlag, Berlin, 2017, pp. 227--235.

\bibitem{DKMS}
A. Dujella, M. Kazalicki, M. Miki\'c and M. Szikszai, \textit{There are infinitely many rational
Diophantine sextuples}, Int. Math. Res. Not. IMRN {\bf 2017 (2)} (2017), 490--508.

\bibitem{DKP-sext}
A. Dujella, M. Kazalicki and V. Petri\v{c}evi\'c,
\textit{There are infinitely many rational Diophantine sextuples with square denominators},
J. Number Theory {\bf 205} (2019), 340--346.

\bibitem{DKP-split}
A. Dujella, M. Kazalicki and V. Petri\v{c}evi\'c, \textit{Rational Diophantine sextuples containing two
regular quadruples and one regular quintuple}, Acta Mathematica Spalatensia, to appear.

\bibitem{DP}
A. Dujella and V. Petri\v{c}evi\'{c},
\textit{Strong {D}iophantine triples},
 Experiment. Math. {\bf 17} (2008), 83--89.

\bibitem{Gibbs1}
P. Gibbs,
\textit{Some rational Diophantine sextuples}, Glas. Mat. Ser. III {\bf 41} (2006), 195--203.

\bibitem{GM}
F. Gouv\^{e}a and B. Mazur, \textit{The square-free sieve and the rank of elliptic curves},
J. Amer. Math. Soc. {\bf 4} (1991), 1--23.

\bibitem{HTZ}
B. He, A. Togb\'e and V. Ziegler,
{\it There is no Diophantine quintuple}, Trans. Amer. Math. Soc. {\bf 371} (2019), 6665--6709.

\bibitem{Hea}
T. L. Heath, Diophantus of Alexandria. A Study in the History of Greek Algebra. Powell's Bookstore, Chicago; Martino Publishing, Mansfield Center, 2003.

\bibitem{Knapp}
A. Knapp, Elliptic Curves, Princeton Univ. Press, 1992.

\bibitem{M}
J.-F. Mestre, \textit{Rang de certaines familles de courbes elliptiques d'invariant donn\'e},
C. R. Acad. Sci. Paris {\bf 327} (1998), 763--764.

\bibitem{RS}
K. Rubin and A. Silverberg, \textit{Rank frequencies for quadratic twists of elliptic curves},
Exper. Math. {\bf 10} (2001), 559--569.


\end{thebibliography}
\end{document}